\documentclass[11pt]{amsart}
\usepackage{amsxtra}

\theoremstyle{plain}

\newcommand{\cleqn}{\setcounter{equation}{0}}
\newcommand{\clth}{\setcounter{theorem}{0}}
\newcommand {\sectionnew}[1]{\section{#1}\cleqn\clth}

\newtheorem{theorem}{Theorem}[section]
\newtheorem{lemma}[theorem]{Lemma}
\newtheorem{definition-theorem}[theorem]{Definition-Theorem}
\newtheorem{proposition}[theorem]{Proposition}
\newtheorem{corollary}[theorem]{Corollary}
\newtheorem{definition}[theorem]{Definition}
\newtheorem{example}[theorem]{Example}
\newtheorem{remark}[theorem]{Remark}
\newtheorem{notation}[theorem]{Notation}
\newtheorem{assumption}[theorem]{Assumption}
\newtheorem{lemma-definition}[theorem]{Lemma-Definition}
\newcommand \bth[1] { \begin{theorem}\label{t#1} }
\newcommand \ble[1] { \begin{lemma}\label{l#1} }

\newcommand \bpr[1] { \begin{proposition}\label{p#1} }
\newcommand \bco[1] { \begin{corollary}\label{c#1} }
\newcommand \bde[1] { \begin{definition}\label{d#1}\rm }
\newcommand \bex[1] { \begin{example}\label{e#1}\rm }
\newcommand \bre[1] { \begin{remark}\label{r#1}\rm }

\newcommand \bnota[1] {\begin{notation}\label{n#1}\rm }
\newcommand \bas[1] { \begin{assumption}\label{a#1}\rm }

\newcommand {\eth} { \end{theorem} }
\newcommand {\ele} { \end{lemma} }

\newcommand {\epr} { \end{proposition} }
\newcommand {\eco} { \end{corollary} }
\newcommand {\ede} { \end{definition} }
\newcommand {\eex} { \end{example} }
\newcommand {\ere} { \end{remark} }
\newcommand {\enota} { \end{notation} }
\newcommand {\eas} {\end{assumption}}

\newcommand \thref[1]{Theorem \ref{t#1}}
\newcommand \leref[1]{Lemma \ref{l#1}}
\newcommand \prref[1]{Proposition \ref{p#1}}

\newcommand \deref[1]{Definition \ref{d#1}}
\newcommand \exref[1]{Example \ref{e#1}}
\newcommand \reref[1]{Remark \ref{r#1}}
\newcommand \lb[1]{\label{#1}}

\def \calM {{\mathcal M}}

\def \de {\delta}

\def \de {\delta}

\def \la {\lambda}

\def \ra  {\rightarrow}           

\def \la {\langle}
\def \ra {\rangle}

\def \Ad { {\mathrm{Ad}} }

\def \g  {\mathfrak{g}}   
\def \h  {\mathfrak{h}}

\def \n  {\mathfrak{n}}

\def \b  {\mathfrak{b}}

\def \la{{\langle}}
\def \ra{{\rangle}}

\DeclareMathOperator \Aut { {\mathrm{Aut}} }

\def \hs {\hspace{.2in}}

\def \lara {\la \, \, , \, \ra}

\def \lara {\la  \, , \, \ra}

\def \hs {\hspace{.2in}}
\def \dw {\dot{w}}

\newcommand{\dias}{\begin{flushright} $\diamondsuit$ \ \ \ \ 
                  \end{flushright} }
                                    
\def \th {\theta}

\def \cth {\cdot_\th}
\def \sc {{\scriptscriptstyle C}}

\def \mc {m_{\sc}}

\def \rk {{\rm rk}}

\def \del {\delta}
\def \tMth {\widetilde{\mathcal M}_\th}

\def \calI {\mathcal I}
\input xy
\xyoption{all}
\begin{document}
\setlength{\baselineskip}{1.2\baselineskip}
\title[On a dimension formula for twisted spherical conjugacy classes]
{On a dimension formula for twisted spherical conjugacy classes in semisimple algebraic groups}
\author[Jiang-Hua Lu]{Jiang-Hua Lu}
\address{Jiang-Hua Lu,
Department of Mathematics,
Hong Kong University,
Pokfulam Rd., Hong Kong}
\email{jhlu@maths.hku.hk}
\date{}
\begin{abstract} Let $G$ be a connected semisimple algebraic group over an algebraically closed
field of characteristic zero,
and let $\th$ be an  automorphism of $G$. We give a characterization of 
$\th$-twisted spherical conjugacy classes 
in $G$ by a formula for their dimensions in terms of certain elements in the Weyl group of $G$, 
generalizing a result of
N. Cantarini, G. Carnovale, and M. Costantini when $\th$ is the identity automorphism. 
For $G$ simple and $\th$ an outer automorphism of $G$, we also classify 
the Weyl group elements that appear in the dimension formula.
\end{abstract}
\maketitle

\sectionnew{Introduction}\lb{sec-intro}

\subsection{The main results}\lb{subsec-result}
Let $G$ be a connected semisimple algebraic group over an algebraically closed  field ${\bf k}$ of characteristic zero.
For an automorphism $\th \in {\rm Aut}(G)$ of $G$, define the $\theta$-twisted conjugation of $G$ on itself by
$g_1 \cdot_\theta g = g_1 g \theta(g_1)^{-1}$ for $g_1, g \in G$,
and call its orbits the $\th$-twisted conjugacy classes in $G$. A $\th$-twisted conjugacy class $C$ in $G$ is said to be
spherical if a Borel subgroup of $G$ has an open orbit in $C$.

Fix a Borel subgroup $B$ of $G$ and a maximal torus $H\subset B$, and let ${\rm Aut}'(G) = \{\th \in {\rm Aut}(G): \th(B) = B, \th(H) = H\}$.
Throughout the paper, we assume that $\th \in \Aut'(G)$ 
(see \reref{re-th-any}).

Let $W=N_G(H)/H$ be the Weyl group, where
$N_G(H)$ is the stabilizer subgroup of $H$ in $G$, and let $l$ be the length function on $W$. 
For $w\in W$, denote by $\rk(1-w\th)$ the rank of the linear operator $1-w\th$ on the Lie algebra $\h$ of $H$.
For a $\th$-twisted conjugacy class $C$ in $G$, let
$\mc$ be the unique element in $W$ such that $C \cap (B\mc B)$ is dense in $C$. In the first part of the paper, we prove the following
characterization of $\th$-twisted spherical conjugacy classes in $G$.

\bth{th-main}
For $\th \in \Aut'(G)$, a $\th$-twisted conjugacy class $C \subset G$ is spherical if and only if
\begin{equation}\lb{eq-formula}
\dim C = l(\mc) + \rk(1-\mc \th).
\end{equation}
\eth

When $\th={\rm Id}_G$, the identity automorphism of $G$, \thref{th-main} is 
proved by N. Cantarini, G. Carnovale, and M. Costantini in \cite{CCC} by
a case-by-case checking that depends  on the classification of all
spherical conjugacy classes in $G$ (for $G$ simple).
Formula \eqref{eq-formula} is then used in \cite{CCC} to prove the De Concini-Kac-Procesi conjecture
on representations of the quantized enveloping algebra of $G$ at roots of unity over spherical conjugacy classes. 
A different proof of \thref{th-main} for $\theta = {\rm Id}_G$, which is also valid when 
the characteristic of ${\bf k}$ is an odd good prime for
$G$, is given by G. Carnovale in \cite{Ca:MathZ}, where
the proof does not require a classification of
spherical conjugacy classes in $G$ but it also depends on some case-by-case computations. 
When $\th^2 = {\rm Id}_G$ and $C$ is the $\th$-twisted conjugacy class through the identity element
of $G$, 
\eqref{eq-formula} follows from standard results on 
symmetric spaces (see $\S$\ref{subsec-symmetric-1}).
 
In $\S$\ref{sec-proof}, we give a direct proof of \thref{th-main}.

For $\th = {\rm Id}_G$, the elements $\mc$ play an important role in the study of spherical conjugacy classes. In particular, it is shown by M. Costantini \cite{Co} that
the coordinate ring of a spherical conjugacy class $C$ as a $G$-module is almost entirely
determined by $\mc$ (see \cite[Theorem 3.22]{Co}).
For $G$ simple and of classical type and for $\th = {\rm Id}_G$, the element $\mc$ for every conjugacy 
class in $G$ is computed explicitly in \cite{Chan-thesis}. 
The second part of the paper concerns the set
\begin{equation}\lb{eq-MCth}
\tMth = \{\mc: \; C \; \mbox{is a} \; \mbox{$\th$-twisted conjugacy class in} \; G\} \subset W 
\end{equation}
for an arbitrary $\th \in \Aut'(G)$.
The set $\tMth$ depends only  on the isogeny class of
$G$ (\cite[Remark 2]{CLT}) and the automorphism of the Dynkin diagram of $G$ induced by $\th$ (\reref{re-th-any}). 
Let 
\begin{align}\lb{eq-Mth}
\calM_\th=\{m \in W: \,\; &m \;\mbox{is the unique maximal length element}\\
\nonumber
&\mbox{in its $\th$-twisted conjugacy class in} \; W\}.
\end{align}  
(See $\S$\ref{subsec-Mth}).  By 
\cite[Corollary 2.15]{CLT},
$\tMth \subset \calM_\th.$

For $G$ simple and  $\th$ an inner automorphism of $G$, it is shown in \cite[$\S$3]{CLT} that
$\tMth=\calM_\th$ and 
elements in $\calM_\th$ are classified in \cite[$\S$3]{CLT} using results from \cite{CCC, Ca:MathZ}. For $G$ simple and 
$\th$ an outer automorphism of $G$, we prove in \thref{th-I} that, again, $\tMth=\calM_\th$, and we give 
in \prref{pr-I} the
complete list of elements in $\calM_\th$. It turns out that if $\th$ induces an order $2$ automorphism of the Dynkin diagram,
the list of elements in $\calM_\th$ coincides with that of  T. A. Springer in \cite[Table 2]{S87}, and if $G = D_4$
and $\theta$ has order $3$, $\calM_\th$ has two elements. 
The classification of elements in $\tMth$ gives restrictions on the possible 
dimensions of $\th$-twisted conjugacy classes in $G$. See \exref{ex-D4}.


\subsection{Notation} Let $\Delta_+$ and $\Gamma\subset \Delta_+$ be the sets of positive and simple roots determined by $(B, H)$ and  
write $\alpha > 0$ (resp. $\alpha < 0$) for  $\alpha \in \Delta_+$ (resp. $\alpha \in -\Delta_+$). 
Let $N$ and $N_-$ be respectively the uniradicals of $B$ and the opposite Borel subgroup $B_-$.  The Lie algebras of 
$G, B, H, N$, and $N_-$ are respectively denoted by $\g, \b, \h, \n$, and $\n_-$.
For $\alpha \in \Delta_+$, $s_\alpha$ denotes the corresponding reflection in $W$. For each $w \in W$, we fix a representative $\dw$ of $w$ in  $N_G(H)$. 

For $\th \in \Aut'(G)$, we use the same letter to denote the action of $\th$ on $\Delta_+$,
and when necessary,
we write $\th \in {\rm Aut}(\Gamma)$ to indicate that $\th$ is regarded as an automorphism of the Dynkin diagram.
The induced action of $\th$ on $\g$ is also denoted by $\th$.

For $g \in G$,  $\Ad_g$ denotes both the conjugation on $G$ by $g$ and the induced
map on $\g$.
For a set $V$ and a map $\sigma: V \to V$,  we let $V^\sigma = \{x \in V: \sigma(x) = x\}$.

\bre{re-th-any}
For an arbitrary $\th_1 \in \Aut(G)$,
there exists $g_0 \in G$ such 
that $\Ad_{g_0}(B) = \th_1(B)$ and $\Ad_{g_0}(H) = \th_1(H)$, so  $\th = \Ad_{g_0}^{-1}\circ \th_1 
\in {\rm Aut}'(G)$, and  the right translation by $g_0$ 
in $G$ maps $\th_1$-twisted conjugacy classes in $G$ to $\th$-twisted conjugacy classes in $G$.
We can thus assume throughout the paper that $\th \in \Aut'(G)$. Moreover, if $\th$ and $\th' \in \Aut'(G)$
are in the same inner class, i.e., if they induce the same automorphism on the Dynkin diagram, then 
$\th = \Ad_h \circ \th'$ for some $h \in H$, and it follows that $\tMth = \widetilde{\calM}_{\th'}$.
\ere

\subsection{Acknowledgment} The author would like to thank G. Carnovale, 
K. Y. Chan and  G. Lusztig for discussions and for answering questions. The author is especially grateful to
G. Carnovale and M. Costantini for pointing out some of the results in \cite{Co:bad}.
Research in this paper was partially supported by HKRGC grants 703405 and 703707.

\sectionnew{Proof of \thref{th-main}}\lb{sec-proof}

\subsection{Two lemmas on $B$-orbits in $G$}\lb{subsec-lemma-B}
 Recall that $\cth$ denotes the
$\th$-twisted conjugacy action of $G$ on itself. For $g \in G$, let $B_g$ be the stabilizer subgroup of $B$ at $g$.
The following generalization of 
\cite[Theorem 5]{CCC} is proved in \cite[Theorem 4.1]{Co:bad}. We include the (short) proof for the convenience of the
reader and to make the proof of \thref{th-main} self-contained.

\ble{le-stabi-B}\cite{Co:bad} 
For any $w \in W$ and $g \in wB$, one has $B_g \subset H^{w\th}(N \cap \Ad_{\dw}(N))$.
Consequently,
\[
\dim B \cth g \geq l(w) + \rk (1-w\th).
\]
\ele

\begin{proof}
Let $b =n_1n_2h\in B_g$, where $h \in H$, $n_1 \in N \cap \Ad_{\dw}(N_-)$ and
$n_2 \in N \cap \Ad_{\dw}(N)$. It follows from $bg = g \theta(b)$ and the unique decomposition
$B w B = (N \cap \Ad_{\dw}(N_-)) \dw B$ that
$n_1 = 1$ and $w\th(h) = h$. Thus $B_g \subset H^{w\th}(N \cap \Ad_{\dw}(N))$, and 
\begin{align*}
\dim B \cth g &= \dim B - \dim B_g \geq \dim B - \dim (N \cap \Ad_{\dw}(N)) - \dim H^{w\th}\\
&=l(w) + \rk (1-w\th).
\end{align*}
\end{proof}

\ble{le-bg}
If $w \in W$ and $g \in wB$ are such that $B \cth g$ is open in $G \cth g$, then
$B_g$ is an open subgroup of $H^{w\th}(N \cap \Ad_{\dw}(N))$. 
\ele

\begin{proof} Let $\g_g=\{x \in \g: \Ad_g \th(x) = x\}$ be the stabilizer subalgebra of $\g$ at $g$ 
for the $\th$-twisted conjugation action, and 
let $\b_g = \b \cap \g_g$. By 
\leref{le-stabi-B}, $\b_g \subset \h^{w\th} + \n \cap 
\Ad_{\dw}(\n)$. It remains to prove that $\h^{w\th} + \n \cap \Ad_{\dw}(\n)\subset \b_g$. 

Let $x_0 \in \h^{w\th}$ and $x_+ \in \n \cap \Ad_{\dw}(\n)$,
and let $z = (\Ad_g \th)^{-1}(x_+ + x_0) - (x_+ + x_0)$ so that 
$\Ad_g \th (z + x_+ + x_0) = x_+ + x_0$. Using the fact that $\Ad_b (x_0) -x_0 \in \n$ for any $b \in B$, 
one sees that $z \in \n$.  We now show that 
$z = 0$. To this end, let $\lara$ be the Killing form of $\g$.
Since $B \cth g$ is open in $G \cth g$,
the inclusion $\b \hookrightarrow \g$ induces an isomorphism 
$\b/\b_g \cong \g/\g_g$. Thus for any $y \in \g$, there exists $y' \in \b$ such that
$y-y' \in \g_g$, and, using  $\la z,\; y' \ra = 0$, one has
\begin{align*}
\la z, \, y \ra &=\la z + x_+ + x_0, \; y-y' \ra - \la x_+ + x_0, \; y-y'\ra \\
&=\la z + x_+ + x_0, \; y-y' \ra - \la \Ad_g \th(z+x_+ + x_0), \; \Ad_g\th(y-y')\ra=0.
\end{align*}
It follows that $z = 0$ and hence $x_+ + x_0 \in \b_g$. Therefore $\b_g = \h^{w\th} + \n \cap 
\Ad_{\dw}(\n)$.
\end{proof}

\subsection{Proof of \thref{th-main}}\lb{subsec-proof}
Let $C$ be a $\th$-twisted conjugacy class in $G$. 
Assume first that $\dim C = l(\mc) + \rk(1-\mc \th)$. By \leref{le-stabi-B}, every $B$-orbit in $C \cap (B \mc B)$ is open in
$C$, so $C$ is spherical.

Assume that $C$ is spherical. Let $g \in C$ be such that $B \cth g$ is open  in $C$, and let 
$g \in BwB$ with $w \in W$. Then $C \cap (BwB) \supset B\cth g$ is dense in $C$, so $w = \mc$. By \leref{le-bg},
\[
\dim C = \dim \b -\dim \b_g = l(\mc) + \rk(1-\mc\th).
\]
This finishes the proof of \thref{th-main}.

\bre{re-vc} For $\th = {\rm Id}_G$,  \leref{le-bg} is also proved in \cite{Ca:MathZ} by some case-by-case arguments. On the other hand,
the arguments in \cite{Ca:MathZ} are 
valid 
when the characteristic of ${\bf k}$ is an  odd good prime for $G$, while our proof of \leref{le-bg} is valid when the Killing for of $\g$ is non-degenerate 
 and when one has the identifications of tangent spaces
$T_g (B \cth g) \cong \b/\b_g$ and $T_g (G\cth g) 
\cong \g/\g_g$, which hold when ${\bf k}$ is of characteristic zero.
\ere


\subsection{The case of symmetric spaces}\lb{subsec-symmetric-1}
Assume that $\th \in {\rm Aut}'(G)$ is an involution, and let $K = G^\th$ be the fixed point subgroup of $\th$ in $G$.
Then the $\th$-twisted conjugacy class $C$ of the identity element of $G$ is isomorphic to the symmetric space
$G/K$, and it is well-known \cite{S85} that $G/K$ is spherical. In this case, formula \eqref{eq-formula}
for the dimension of $G/K$ follows from results in \cite{S85}. Indeed, using the notation in \cite[$\S$5]{S85}, let
$v^o$ be the unique open $B$-orbit in $G/K$ and let $w^o = \phi(v^o) \in W$. Then $w^o = \mc$, and it is easy to see from
\cite[Corollary 4.9]{S85} that $\dim G/K = \frac{1}{2} {\rm Card}(C_{v^o}^{\prime\prime}) + {\rm Card}(I_{v^o}^n) + l(w^o) + \rk (1-w^o\th)$,
where the notation is as on \cite[Page 535]{S85}. By \cite[Theorem 5.2(i)]{S85}, $C_{v^o}^{\prime\prime} \cap \Gamma = \emptyset$.
For every $\beta > 0$, writing 
$\beta =\beta_1 + \beta_2$, where $\beta_1$ is in the linear span of $\Pi\subset \Gamma$ 
in the notation of \cite[Theorem 5.2(ii)]{S85}
and $\beta_2$ is in the linear span of $\Gamma \backslash \Pi$, one has $w^o\th(\beta) = \beta_1 + w^o\th(\beta_2)$, so 
by \cite[Theorem 5.2(ii)]{S85}, $w^o\th(\beta) > 0$ implies that $\beta_2 = 0$ and thus $w^o\th(\beta) =\beta$. This shows that
$C_{v^o}^{\prime\prime}=\emptyset$ and that every $\beta \in I_{v^o}$ is in the linear span of 
$\Pi$, which, by \cite[Theorem 5.2(i)]{S85}, consists of all simple compact imaginary roots. It follows that $I_{v^o}^n=\emptyset$.
Thus $\dim G/K = l(w^o) + \rk (1-w^o\th)$.

\sectionnew{The elements $\mc$}\lb{sec-mc}

\subsection{Properties of $m \in \calM_\th$}\lb{subsec-Mth}
Any $\delta \in {\rm Aut}(\Gamma)$ induces an automorphism on the Weyl group $W$, also denoted by $\delta$, by 
$\del(w) = \de \circ w \circ \del^{-1}:
\h \to \h$. Define the $\del$-twisted conjugacy of $W$ on itself by $w \cdot_\delta v = w v \del(w)^{-1}$ for $w, v \in W$
and call its orbits $\del$-twisted conjugacy classes in $W$. 
Let $w_0$ be the longest element in $W$, and let $\delta_0$ be the automorphism of 
$\Gamma$ given by $\delta_0(\alpha) = -\alpha$ for $\alpha \in \Gamma$. The automorphism on $W$ induced by
$\delta_0$ is then given by  $\delta_0(w) = w_0 w w_0$ for $w \in W$. 

Throughout this section, $\th \in {\rm Aut}(\Gamma)$, and 
 $\calM_\th \subset W$ is defined as in \eqref{eq-Mth}.

\ble{le-m-0}
If $m \in \calM_\th$, then $\th(m) = \delta_0(m) = m$.
\ele

\begin{proof}
Let $m \in \calM_\th$. Then $\th(m) =m^{-1} m \th(m)$ is in the same $\th$-twisted conjugacy class 
as $m$, and $l(\th(m)) = l(m)$. Thus $\th(m) = m$. Similarly,  
since $\th$ permutes the simple roots, 
$\theta(w_0) = w_0$. It follows that $w_0 m w_0$ and $m$ are in the same $\th$-twisted conjugacy
class in $W$. Since $l(w_0 m w_0) = l(m)$, one has
$w_0 m w_0 = m$.
\end{proof}

For $I \subset\Gamma$, let  $w_{0, I}$ be the
longest element in the subgroup $W_I$ of $W$ generated by $I$. 

The following \leref{le-m} is proved in \cite[$\S$3]{CLT}  when
$\th$ is the identity automorphism of $\Gamma$.

\ble{le-m}
If $m \in \calM_\th$, then
$w_0m =mw_0 = w_{0, I}$, where $I = \{\alpha \in \Gamma: m\th(\alpha) = \alpha\}$. In particular,
$I$  is both $\delta_0$ and $\th$ invariant, 
and $\delta_0\th(\alpha) = -w_{0, I}(\alpha)$ for every
$\alpha \in I$. 
\ele

\begin{proof} 
Let $\delta = \delta_0 \th \in {\rm Aut}(\Gamma)$. Then the map $W \to W: w \mapsto ww_0$ maps
$\th$-twisted conjugacy classes in $W$ to $\delta$-twisted conjugacy classes in $W$. 

Let $m \in \calM_\th$, 
and let $x =mw_0$. Then $x$ 
is a unique minimal length 
element in its $\delta$-twisted conjugacy class  in $W$. Let
$x = s_{\alpha_1} s_{\alpha_2} \cdots s_{\alpha_k}$ be a reduced word for $x$.
 Let $I' = \{\alpha_1, \alpha_2, \ldots, \alpha_k\}$. Then $x \in W_{I'}$. We first show
that $x = w_{0, I'}.$ To this end, 
 it is enough to show that $x(\alpha_j) <0$ for every $1 \leq j \leq k$. 
Since $x s_{\alpha_k} < x$, we already know that $x(\alpha_k) <0$. If $k = 1$, we are done.
Suppose that $k \geq 2$. Let $\beta_k = \delta^{-1}(\alpha_k) \in \Gamma$,  and let 
\begin{equation}\lb{eq-x1}
x_1 =
s_{\beta_k} x \delta(s_{\beta_k})=s_{\beta_k} x s_{\alpha_k} = s_{\beta_k} s_{\alpha_1}  \cdots s_{\alpha_{k-1}}.
\end{equation}
Since $k$ is the minimal length of elements in the $\delta$-twisted conjugacy class of $x$ in $W$, we have 
$l(x_1) \geq k$. It follows from
\eqref{eq-x1} that $l(x_1) \leq k$, so $l(x_1) = k$. Since $x$ is the unique element in its $\delta$-twisted conjugacy class
in $W$ with length $k$, we have $x_1 = x$.  
In particular, 
$x = x_1 =s_{\beta_k} s_{\alpha_1}  \cdots s_{\alpha_{k-1}}$ is a reduced word for $x$, so
$x(\alpha_{k-1}) <0$. Repeating this process, we see that $x(\alpha_j)< 0$
for every $1 \leq j \leq k$. Thus $x=w_0m=mw_0 = w_{0, I'}$. It follows from
\leref{le-m-0} that $\delta_0(I') = \th(I') = I'$.

We now  show that $I' = I$. For any $\alpha \in I'$, since $m(\alpha)=w_0w_{0, I'}(\alpha)> 0$, one has
$l(\th^{-1}(s_\alpha) m s_\alpha) \geq l(m)$. Since $m \in \calM_\th$, one has $\th^{-1}(s_\alpha) m s_\alpha = m$, 
so $\th^{-1}(\alpha) = m(\alpha)$ and $\alpha \in I$. Conversely, let $\alpha \in I$. If $\alpha \notin I'$, then 
$w_0m(\alpha) = w_{0, I'}(\alpha) > 0$, so $m(\alpha) < 0$, contradicting the fact that $m(\alpha) = \th^{-1}(\alpha) > 0$.
Thus $I' = I$. It follows from the definition of $I$ that
$\delta_0\th(\alpha) = -w_{0, I}(\alpha)$ for every
$\alpha \in I$.
\end{proof} 

An element $w \in W$ is said to be a $\th$-twisted involution if $\th(w) = w^{-1}$.

\bco{co-m2}
Every $m \in \calM_\th$ is both an involution and a $\th$-twisted involution.
\eco

\begin{proof}
Let $m \in \calM_\th$ and let the notation be as in \leref{le-m}. Then $m^2 = w_0w_{0, I} w_{0, I} w_0 = 1$.
Since $\th(m) = m$, one also has $\th(m) = m^{-1}$.
\end{proof}

\bde{de-Property1}
A subset  $I$ of $\Gamma$ is said to have Property (1) if $I$ is both $\delta_0$ and $\th$ invariant and if
$\delta_0\th(\alpha) = -w_{0, I}(\alpha)$ for all $\alpha \in I$.
\ede

By  \leref{le-m}, every $m \in \calM_\th$ is of the form $m = w_0 w_{0, I}$ for some $I \subset \Gamma$ with
Property (1).
The following \deref{de-Property2} is inspired by \cite[Lemma 4.1]{Ca:MathZ}.


\bde{de-Property2} For a subset $I$ of $\Gamma$, an $\alpha \in I$ is said to be isolated if $\la \alpha, \alpha'\ra = 0$
for every $\alpha' \in I\backslash\{\alpha\}$.
A subset  $I$ of $\Gamma$ is said to have Property (2) if for every isolated $\alpha \in I$, there is no
$\beta \in \Gamma\backslash\{\alpha\}$ with the following properties

(a) $\la \alpha, \alpha \ra = \la \beta, \beta\ra$ and $\la \beta, \alpha\ra \neq 0$;

(b) $\la \beta, \alpha'\ra = 0$ for all $\alpha' \in I\backslash \{\alpha\}$;

(c) $\delta_0\th(\beta) =\beta$.
\ede

\ble{le-Property2}
For every $m \in \calM_\th$,  $I_m = \{\alpha \in \Gamma: m\th(\alpha) = \alpha\}\subset \Gamma$ has Property (2).
\ele

\begin{proof} 
Let $m \in \calM_\th$. Suppose that
$\alpha \in I_m$ is isolated and that there exists $\beta \in \Gamma\backslash\{\alpha\}$ with properties (a), (b), and (c)
in \deref{de-Property2}. Let $I_m'=I_m\backslash\{\alpha\}$. Since $\alpha \in I_m$ is isolated, $w_{0, I_m} = s_\alpha w_{0, I_m'}$, 
so by (b) and (c), $m\th(\beta) = w_{0, I_m} w_0\th(\beta) = -s_\alpha w_{0, I_m'}(\beta) = -s_\alpha(\beta)$, and 
\[
s_\alpha s_\beta m s_{\th(\beta)}s_{\th(\alpha)} =s_\alpha s_\beta s_{m\th(\beta)} m s_{\th(\alpha)} = 
s_\alpha s_\beta s_\alpha s_\beta s_\alpha m s_{\th(\alpha)}.
\]
By (a), $s_\alpha s_\beta s_\alpha s_\beta s_\alpha=s_\beta$, so $s_\alpha s_\beta m s_{\th(\beta)}s_{\th(\alpha)}
=s_\beta m s_{\th(\alpha)} = s_\beta s_\alpha m$. Since  $m^{-1}(\alpha) = \th(\alpha) > 0$, $l(s_\beta s_\alpha m) \geq l(m)$.
Since $s_\alpha s_\beta m$ is in the same $\th$-twisted conjugacy class as $m$,  we have $s_\beta s_\alpha m=m$,
or $s_\alpha s_\beta = 1$, which is a contradiction.
\end{proof}

\subsection{The classification of $m \in \calM_\th$}\lb{subsec-classification}

For $\th \in {\rm Aut}(\Gamma)$, let $\calI_\th$ be the collection of all subsets $I$ of $\Gamma$ that have Properties (1) and (2). Note that the empty set $\emptyset$ is
always in $\calI_\th$. 
Also note that if $\th \in {\rm Aut}(\Gamma)$ is not the identity automorphism, then 
$\Gamma$ does not have Property (1),
so
$\Gamma \notin \calI_\th$.

\bpr{pr-I} 1) For $G = D_4$ and $\th \in {\rm Aut}(\Gamma)$ of order $3$,
$I \in \calI_\th$ if and only if $I = \emptyset$ or $I = \{\alpha_2\}$, where 
$\alpha_2$ is the simple root that is not orthogonal to
any of the other three.

2) For $G$ simple and $\th \in {\rm Aut}(\Gamma)$ of order $2$,
the list for $I \in \calI_\th$ is the same as that given in \cite[Table 2]{S87}, namely, 
either $I$ is the empty set or $I$ is
one the following:

$A_{2n}, \;n \geq 1,\;\th =\delta_0$:  no non-empty $I$ in $\calI_\th$.

$A_{2n+1}, \; n \geq 1, \;\th =\delta_0$: $I = \{\alpha_{2l+1}: 0 \leq l \leq n\}$.

$D_4: \; I =\{\alpha_2\} \cup \Gamma(2, \th)$, where $\Gamma(2, \th)$ is the $\th$-orbit in $\Gamma$ with $2$ elements.

$D_{2n}, \;n > 2, \;\th(\alpha_{2n-1}) =\alpha_{2n}$: $I_l = \Gamma\backslash\{\alpha_1, \alpha_2, \ldots, \alpha_{2l-1}\}$
for $1 \leq l \leq n-1$.

$D_{2n+1}, \;n \geq 2, \;\th = \delta_0$: $I_l = \Gamma\backslash\{\alpha_1, \alpha_2, \ldots, \alpha_{2l-1}\}$
for $1 \leq l \leq n$.

$E_6, \; \th = \delta_0$: $I = \{\alpha_3, \alpha_4, \alpha_5\}$ with the simple roots labeled as 
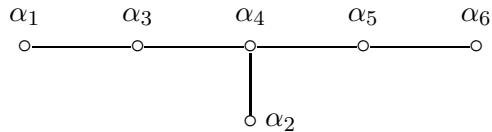
\[
\xy
(0,4)*{\alpha_1};(15,4)*{\alpha_3};(30,4)*{\alpha_4};(45,4)*{\alpha_5};(60,4)*{\alpha_6};(34,-10)*{\alpha_2};
(0,0)*{\circ}="1";(15,0)*{\circ}="3";(30,0)*{\circ}="4";(45,0)*{\circ}="5";(60,0)*{\circ}="6";(30,-10)*{\circ}="2";
{\ar@{-} "1"; "3"  }; 
{\ar@{-} "2"; "4"  }; 
{\ar@{-} "3"; "4"  }; 
{\ar@{-} "4"; "5"  }; 
{\ar@{-} "5"; "6"  }; 
\endxy
\]
\epr

\begin{proof} 1) is easy to deduce and 2) is proved case-by-case. We omit the details. 
\end{proof}
 
By \leref{le-m} and \leref{le-Property2}, we have the well-defined map 
\[
\psi: \;\;\;\; \calM_\th \longrightarrow  \calI_\th: \;\;\; m \longmapsto I_m =\{\alpha \in \Gamma: m\th(\alpha) = \alpha\}.
\]
Since $m = w_0 w_{0, I_m}$ for every $m \in \calM_\th$, the map $\psi$ is injective. 

Assume now $\th \in {\rm Aut}'(G)$, and let
$\tMth \subset W$ be defined as in \eqref{eq-MCth}.
By \reref{re-th-any},
$\tMth$ depends only on the corresponding  $\th \in \Aut(\Gamma)$.
Let $\widetilde{\psi}: \tMth \to \calI_\th$ be the restriction of $\psi$ to $\tMth \subset \calM_\th$.

\bth{th-I}
For $G$ simple and $\th\in \Aut'(G)$ an outer automorphism of $G$, the map
$\widetilde{\psi}: \tMth \to \calI_\th$ is  bijective. Consequently,
\[
\tMth = \calM_\th = \{w_0w_{0, I}: I \in \calI_\th\}.
\]
\eth

\begin{proof} It is enough to prove that $\tilde{\psi}$ is surjective, and we may assume that $G$ is adjoint. 

First assume that $\th \in {\rm Aut}(\Gamma)$ has order $2$, and let $I \in \calI_\th$. By \prref{pr-I}, $I$ is in \cite[Table 2]{S87},
so
$(I, \delta_0\th)$ is {\it admissible} in the sense of
\cite[No. 2.2]{S87}. By \cite[No. 4 and No. 5]{S87}, there exists $h \in H$ such that 
$\Ad_h \th \in {\rm Aut}(G)$ is an involution and that  $w_0w_{0, I} = \mc$, where
$C$ is the $\th$-twisted conjugacy class through $h$.
In particular, 
$w_0 w_{0, I} \in \tMth$.

It remains to consider the case of $G=D_4$ with $\th \in {\rm Aut}(\Gamma)$ having order $3$. 
It is clear that $w_0 = \mc$ if $C$ is the $\th$-twisted conjugacy class of $\dw_0$, so $w_0 \in \tMth$.
We only need to show that $w_0s_2 \in \tMth$. To this end, we may, by \reref{re-th-any}, assume that
$\th \in \Aut'(G)$ is a diagram automorphism of $G$ in the sense that $\th \circ x_\alpha = x_{\th \alpha}$ for $\alpha \in \Gamma$, where
for each $\alpha \in \Gamma$, $x_\alpha: {\bf k}_a \to G$ is a fixed choice of 
one-parameter root subgroup corresponding to $\alpha$.
In particular, $\th^3 = {\rm Id}_G$. Let $C_e$ be the $\th$-twisted conjugacy class through the identity element $e$ of $G$.
It is well-known that $\g^\th$ is of type $G_2$  \cite[Chapter 24]{port} so it is $14$-dimensional. Thus
\[
\dim C_e =\dim G - 14 =14
=l(w_0s_2) + \rk(1-w_0s_2\th).
\]
Since $l(w_0) + \rk(1-w_0\th) = 16$, we know by \leref{le-stabi-B} that
$m_{\sc_e} \neq w_0$ so 
$m_{\sc_e} = w_0s_2$. In particular, $w_0s_2 \in \tMth$ and $C_e$ is spherical.
See \cite[$\S$4.5]{Co:bad} for another proof of the fact that $w_0s_2 \in \tMth$ and that 
$C_e$ is spherical. 
\end{proof}

\bex{ex-D4}
Let $G = D_4$ be of adjoint type, and let $\th \in \Aut'(G)$ be a triality automorphism of $G$ as in the proof of 
\thref{th-I}. Since $l(w_0s_2) + \rk (1-w_0s_2\th) = 14$ and $l(w_0) + \rk (1-w_0\th)=16$, 
$\dim C \geq 14$ for every $\th$-twisted conjugacy class $C$ 
in $G$, and, by \thref{th-main}, $\dim C=14$ or $16$ if $C$ is spherical.

Recall from \cite{S:twisted} that  a $\th$-twisted conjugacy class is semisimple  if it contains an element in $H$.
For $h \in H$, let $C_h \subset G$ be the $\th$-twisted conjugacy class of $h$. 
Label the simple roots as $\Gamma =\{\alpha_j: 1 \leq j \leq 4\}$ such that
$\th(\alpha_2) = \alpha_2$,  $\th(\alpha_1) = \alpha_3, \th(\alpha_3) = \alpha_4$, and $\th(\alpha_4) = \alpha_1$. We now show
that if $h^{\alpha_2} =h^{\alpha_1}h^{\alpha_3}h^{\alpha_4} = 1$, then 
$m_{\sc_h} = w_0s_2$ and $C_h$ is spherical, and otherwise, $m_{\sc_h} = w_0$ and $\dim C_h \geq 20$, so
$C_h$ is 
not spherical. Here, for a character $\mu$ on $H$, $h^\mu$ denotes the 
value of $\mu$ on $h$.

Label the positive roots in $\Delta_+\backslash \Gamma$ as
\begin{align*}
&\alpha_5 =\alpha_1+\alpha_2, \hs \alpha_6 = \alpha_2 + \alpha_3, \hs \alpha_{7} = \alpha_2 + \alpha_4,\\
&\alpha_8 = \alpha_1 + \alpha_2 + \alpha_3, \hs \alpha_9 =\alpha_2 + \alpha_3 + \alpha_4,\hs
\alpha_{10} = \alpha_1 + \alpha_2 + \alpha_4,\\
& \alpha_{11} = \alpha_1 + \alpha_2 + \alpha_3 + \alpha_4, \hs \alpha_{12}
 = \alpha_1 + 2\alpha_2 + \alpha_3 + \alpha_4.
\end{align*}
Then $\{\alpha_1, \alpha_3, \alpha_4\}, \{\alpha_5, \alpha_6, \alpha_7\}$ and $\{\alpha_8, \alpha_9, \alpha_{10}\}$ are
the three $\th$-orbits in $\Delta_+$ of size $3$ and $\th(\alpha_{11}) = \alpha_{11}$ and $\th(\alpha_{12}) = \alpha_{12}$.
Note that the sets $\{\alpha_1, \alpha_3, \alpha_4, \alpha_{12}\}$, 
$\{\alpha_5, \alpha_6, \alpha_7, \alpha_{11}\}$, and $\{\alpha_8, \alpha_9, \alpha_{10}, \alpha_{2}\}$ consist of 
strongly orthogonal roots, and, with $s_j$ denoting the reflection in $W$ defined by $\alpha_j$ for $1 \leq j \leq 12$,
$w_0 = s_1s_3s_4s_{12} = s_5s_6s_7s_{11} = s_8s_9s_{10}s_{2}.$

Recall that the stabilizer subalgebra of $\g$ at $h$ is $\g_h = \g^{\Ad_g\th}$. Since $\dim \h^{\Ad_h\th} = \dim \h^\th = 2$, one has 
$\dim \g_h = 2 + 2 n$, where $n = \#\{i \in \{1, 2, 5, 8, 11, 12\}: \lambda_i(h) = 1\}$, with
$\lambda_i(h) = h^{\alpha_i + \th(\alpha_i) + \th^2(\alpha_i)}$ for $i \in \{1, 5, 8\}$ and
$\lambda_i(h) = h^{\alpha_i}$ for $i \in \{2, 11, 12\}$.
Let $\Lambda(h) = \{\lambda_i(h): i \in \{1, 2, 5, 8, 11, 12\}\}$. Then $\lambda_1(h) = h^{\alpha_1}h^{\alpha_3}h^{\alpha_4}$, $\lambda_2(h) = h^{\alpha_2}$, and 
\[
\Lambda(h) = \{\lambda_1(h), \; \lambda_2(h), \; 
 \lambda_1(h)(\lambda_2(h))^3, \;  (\lambda_1(h))^2(\lambda_2(h))^3, \; \lambda_1(h)\lambda_2(h),
\;  \lambda_1(h)(\lambda_2(h))^2\}.
\]

{\it Case 1:} $h^{\alpha_2} = h^{\alpha_1}h^{\alpha_3}h^{\alpha_4}=1$. In this case, $n = 6$,
$\dim \g_h = 14$, and $\dim C_h = 28-14 = 14$. It follows from \leref{le-stabi-B} that 
$C_h$ is spherical and $m_{\sc_h} = w_0s_2$. Note that in this case,
 $(\Ad_h \th)^3 = {\rm Id}_G$, so
$\g_h$ is again of tyep $G_2$.

{\it Case 2:}  $h^{\alpha_2} \neq 1$ or  $h^{\alpha_1}h^{\alpha_3}h^{\alpha_4}\neq 1$. In this case, $n \leq 5$. In fact, it is easy to see that one can not have $n = 5$ nor $n = 4$, so $n \leq 3$, and $\dim C_h = 28 - \dim \g_h \geq 20$. Thus $C_h$ is not spherical.
We use the approach in \cite[$\S$4.5]{Co:bad} to prove that $m_{\sc_h} = w_0$. First assume that
$h^{\alpha_2} \neq 1$. 
Fix a one-parameter root subgroup $x_\alpha: {\bf k}_a \to G$ for $\alpha \in -\{\alpha_8, \alpha_9, \alpha_{10}\}$ such that
$\th \circ x_{\alpha} = x_{\th(\alpha)}$ for every $\alpha \in -\{\alpha_8, \alpha_9, \alpha_{10}\}$ (recall that
$\th^3 = {\rm Id}_G$).
For $a, b, c, d \in {\bf k}\backslash\{0\}$, let $g = 
x_{-\alpha_2}(a) x_{-\alpha_8}(b) x_{-\alpha_9}(c) x_{-\alpha_{10}}(d) \in G$. Then 
\[
gh\th(g)^{-1} = x_{-\alpha_2}(a-h^{-\alpha_2}a) x_{-\alpha_8}(b-h^{-\alpha_8}d) x_{-\alpha_9}(c-h^{-\alpha_9}b)x_{-\alpha_{10}}
(d-h^{-\alpha_{10}}c).
\]
Choosing $a, b, c, d$ such that $a \neq 0, b-h^{-\alpha_8}d \neq 0, c-h^{-\alpha_9}b \neq 0$ and $d-h^{-\alpha_{10}}c \neq 0$, 
one has
$g h \th(g)^{-1} \in C_h \cap (Bw_0B) \cap B_-$, so
$m_{\sc_h} = w_0$. 
If $h^{\alpha_2} =1$, then $h^{\alpha_1}h^{\alpha_3}h^{\alpha_4}\neq 1$. In this case, $h^{\alpha_{11}} = h^{\alpha_{12}} 
\neq 1$. Using the fact $w_0=s_5s_6s_7s_{11}$ or the fact $w_0 = s_1s_3s_4s_{12}$ and 
arguments similar to the above, one
sees that 
$m_{\sc_h} = w_0$. 
\dias
\eex


\begin{thebibliography}{99}



\bibitem{CCC} N. Cantarini, G. Carnovale, and M. Costantini, {\it Spherical orbits and
representations of ${\mathcal U}_{\epsilon}(\g)$}, Trans. Groups {\bf 10} (1) (2005), 29 - 62.

\bibitem{Ca:MathZ} G. Carnovale, {\it Spherical conjugacy classes and involutions in the Weyl group,} Math. Z. {\bf 260} (1) (2008), 1 - 23.

\bibitem{CLT}
K. Y. Chan, J-H. Lu, and S. To, {\it On intersections of conjugacy classes and Bruhat cells}, to Trans. Groups.

\bibitem{Chan-thesis} K. Y. Chan, {\it MPhil thesis in Mathematics}, The University of Hong Kong, 2010.



\bibitem{Co} M. Costantini, {\it On the coordinate ring of spherical conjugacy classes}, 
 Math. Z. {\bf 264} (2010), 327 - 359.

\bibitem{Co:bad}
M. Costantini, {\it A classification of unipotent conjugacy classes in bad characteristics,} arXiv:0906.5063.




\bibitem{port}
I. Porteous, {\it Clifford algebras and the classical
groups}, Cambridge University Press, 1995.

\bibitem{S85}
T. A. Springer, {\it Some results on algebraic groups with involutions,} 
{\em Algebraic groups and related topics}, Adv. Studies in Pure Math.
{\bf 6} (1985) 525 - 543.

\bibitem{S87}
T. A. Springer, {\it The classification of involutions of simple algebraic groups,} J. Fac. Sci. Univ. Tokyo, Sect. IA, Math
{\bf 34} (1987), 655 - 670.



\bibitem{S:twisted}
T. A. Springer, {\it Twisted conjugacy in simply connected groups},  Trans. Groups. {\bf 11}(3)
(2006), 539 - 545.

\end{thebibliography}
\end{document}